\documentclass[12pt]{amsart}
\usepackage[margin=1in]{geometry}
\usepackage{amsthm,amssymb}
\usepackage{amsmath,xypic}
\usepackage[colorlinks]{hyperref}
\usepackage{epigraph}
\usepackage{graphicx}
\usepackage{fancyhdr} \usepackage{enumitem}
\usepackage{mathrsfs}
\usepackage{natbib}
\let\cite=\citep
\usepackage{zref-clever}
\let\Cref=\zcref

\newcommand{\be}{\begin{equation}}
\newcommand{\ee}{\end{equation}}
\newcommand{\bes}{\begin{equation*}}
\newcommand{\ees}{\end{equation*}}
\newcommand{\bea}{\begin{eqnarray}}
\newcommand{\eea}{\end{eqnarray}}
\newcommand{\beas}{\begin{eqnarray}}
\newcommand{\eeas}{\end{eqnarray}}
\newcommand{\ben}{\begin{note}}
\newcommand{\een}{\end{note}}
\newcommand{\bexl}{\vskip0.1em\noindent\hrulefill\vskip1em\begin{ExerciseList}}
\newcommand{\eexl}{\end{ExerciseList}\hrulefill}

\newcommand{\bthm}{\begin{theorem}}
\newcommand{\ethm}{\end{theorem}}
\newcommand{\bpro}{\begin{prop}}
\newcommand{\epro}{\end{prop}}
\newcommand{\bcor}{\begin{corollary}}
\newcommand{\ecor}{\end{corollary}}
\newcommand{\bcon}{\begin{conjecture}}
\newcommand{\econ}{\end{conjecture}}
\newcommand{\bp}{\begin{proof}}
\newcommand{\ep}{\end{proof}}
\newcommand{\blem}{\begin{lemma}}
\newcommand{\elem}{\end{lemma}}
\newcommand{\bn}{\begin{note}}
\newcommand{\en}{\end{note}}
\newcommand{\benum}{\begin{enumerate}}
\newcommand{\eenum}{\end{enumerate}}
\newcommand{\bed}{\begin{defn}}
\newcommand{\eed}{\end{defn}}
\newcommand{\brem}{\begin{remark}}
\newcommand{\erem}{\end{remark}}

\newcommand{\btik}{\begin{tikzpicture}\begin{axis}[scale=0.5,axis y line=center, axis x line=middle]}
\newcommand{\etik}{\end{axis}\end{tikzpicture}}

\let\mapsto=\longmapsto

\newcommand{\upperRomannumeral}[1]{\uppercase\expandafter{\romannumeral#1}}

\newtheorem{theorem}[equation]{Theorem}      \newtheorem{lemma}[equation]{Lemma}          \newtheorem{corollary}[equation]{Corollary}  \newtheorem{proposition}[equation]{Proposition}

\theoremstyle{definition}

\theoremstyle{definition}
\newtheorem{defn}[equation]{Definition}
\theoremstyle{remark}

\theoremstyle{definition}
\newtheorem{remark}[equation]{Remark}

\numberwithin{equation}{section}

\let\isom=\simeq

\newcommand{\A}{\mathcal{A}}

\newcommand{\C}{{\mathbb C}}

\newcommand{\Pic}{{\rm Pic\,}}

\newcommand{\spec}{{\rm Spec}}

\newcommand{\Z}{{\mathbb Z}}

\renewcommand{\int}{\operatorname{int}}

\renewcommand{\P}{{\mathbb P}}

\newcommand{\alb}{{\rm Alb}}

\newcommand{\mapright}[1]{{\xymatrix{{}\ar[r]^{#1}&{}}}}

\renewcommand{\bpro}{\begin{proposition}}
\renewcommand{\epro}{\end{proposition}}

\zcsetup{cap} \AddToHook{env/theorem/begin}{\zcsetup{countertype={equation=theorem}}}
\AddToHook{env/corollary/begin}{\zcsetup{countertype={equation=corollary}}}
\AddToHook{env/lemma/begin}{\zcsetup{countertype={equation=lemma}}}
\AddToHook{env/proposition/begin}{\zcsetup{countertype={equation=proposition}}}
\AddToHook{env/prop/begin}{\zcsetup{countertype={equation=prop}}}
\AddToHook{env/theoremdef/begin}{\zcsetup{countertype={equation=theoremdef}}}
\AddToHook{env/defn/begin}{\zcsetup{countertype={equation=defn}}}
\AddToHook{env/definition/begin}{\zcsetup{countertype={equation=definition}}}
\AddToHook{env/remark/begin}{\zcsetup{countertype={equation=remark}}}
\AddToHook{env/rem/begin}{\zcsetup{countertype={equation=rem}}}
\AddToHook{env/conj/begin}{\zcsetup{countertype={equation=conj}}}
\AddToHook{env/conjecture/begin}{\zcsetup{countertype={equation=conjecture}}}
\AddToHook{env/question/begin}{\zcsetup{countertype={equation=question}}}
\AddToHook{env/example/begin}{\zcsetup{countertype={equation=example}}}

\zcRefTypeSetup{equation}{
	Name-sg = eq. ,
	name-sg = eq. ,
	Name-pl = eqns. ,
	name-pl = eqns. ,
}

\zcRefTypeSetup{theorem}{
	Name-sg = Theorem ,
	name-sg = theorem ,
	Name-pl = Theorems ,
	name-pl = theorems ,
}

\zcRefTypeSetup{theoremdef}{
	Name-sg = Theorem-Definition ,
	name-sg = theorem-definition ,
	Name-pl = Theorem-Definitions ,
	name-pl = theorem-definitions ,
}

\zcRefTypeSetup{corollary}{
	Name-sg = Corollary ,
	name-sg = corollary ,
	Name-pl = Corollaries ,
	name-pl = corollaries ,
}

\zcRefTypeSetup{proposition}{
	Name-sg = Proposition ,
	name-sg = proposition ,
	Name-pl = Propositions ,
	name-pl = propositions ,
}

\zcRefTypeSetup{prop}{
	Name-sg = Proposition ,
	name-sg = proposition ,
	Name-pl = Propositions ,
	name-pl = propositions ,
}

\zcRefTypeSetup{lemma}{
	Name-sg = Lemma ,
	name-sg = lemma ,
	Name-pl = Lemmas ,
	name-pl = lemmas ,
}

\zcRefTypeSetup{remark}{
	Name-sg = Remark ,
	name-sg = remark ,
	Name-pl = Remarks ,
	name-pl = remarks ,
}

\zcRefTypeSetup{question}{
	Name-sg = Question ,
	name-sg = question ,
	Name-pl = Questions ,
	name-pl = questions ,
}

\zcRefTypeSetup{example}{
	Name-sg = Example ,
	name-sg = example ,
	Name-pl = Examples ,
	name-pl = examples ,
}

\zcRefTypeSetup{defn}{
	Name-sg = Definition ,
	name-sg = definition ,
	Name-pl = Definitions ,
	name-pl = definitions ,
}

\zcRefTypeSetup{conj}{
	Name-sg = Conjecture ,
	name-sg = conjecture ,
	Name-pl = Conjectures ,
	name-pl = conjectures ,
}

\zcRefTypeSetup{conjecture}{
	Name-sg = Conjecture ,
	name-sg = conjecture ,
	Name-pl = Conjectures ,
	name-pl = conjectures ,
}

\begin{document}

\title[]{On the Grothendieck ring of varieties in positive characteristic}\author{Kirti Joshi}\address{Math. department, University of Arizona, 617 N Santa Rita, Tucson
85721-0089, USA.}

\thanks{}\subjclass{}\keywords{}

\begin{abstract}
 This paper proves two theorems
(1) Let $k$ be an algebraically closed field of characteristic $p>0$. I prove (\Cref{th:main}) that if, $p>13$ or $p=11$, then the  isomorphism class of any supersingular elliptic curve is a zero divisor in the ring of smooth, complete $k$-varieties and Bittner relations. In particular, this ring contains zero divisors. The proof proceeds via establishing (in \Cref{th:alb3}) that the Albanese variety functor is a motivic measure. (2) I  prove (\Cref{th:pi1-is-a-motivic-measure}) that the \'etale fundamental group of a smooth, proper variety (in any characteristic) also provides a motivic measure on this ring. In particular, the \'etale fundamental group is a motivic measure on the Grothendieck ring of varieties over complex numbers.
\end{abstract}
\maketitle
\renewcommand{\L}{{\mathbb{L}}}
\renewcommand{\A}{\mathbb{A}}

\tableofcontents
\section{Introduction}

\let\citep=\cite
\newcommand{\gR}{{\rm K}^0(\mathcal{V}_k)}
\newcommand{\gRs}{{\rm K}^0_{\rm bl}(\mathcal{CV}^{sm}_k)}
\newcommand{\B}[2]{{\rm Bl}_{#1}({#2})}
\newcommand{\chr}[1]{{\rm char}(#1)}
\newcommand{\bw}{\bar{w}}

Let $k$ be  an algebraically closed field of characteristic $p>0$. Let $\gR$ be the Grothendieck ring of varieties over $k$. This is the following: for an irreducible variety $X/k$, write $[X]$ for its $k$-isomorphism class and consider the free abelian group generated by isomorphism classes $[X]$ of $k$-varieties with multiplication given product of varieties: $[X_1][X_2]=[X_1\times X_2]$ of $k$-varieties $X_1,X_2$. This construction provides a ring in which $[\emptyset]=0$ and $[\spec(k)]=1$ is the multiplicative identity element. Let $\gR$ be the quotient by the ideal generated by elements of the form $[X-Y]-[X]-[Y]$ where $X\supseteq Y$ is a closed $k$-subvariety of $X$. This makes $\gR$ into a commutative ring with a unit (with $[\spec(k)]=1$)  see \citep{larsen03,bittner01} for the basic properties of this ring.

While the simplicity of  definition is one of the most attractive features of $\gR$, it should be noted that this ring remembers far less algebraic geometry than what one learns  in a first course on algebraic geometry. So one can ask what parts of algebraic geometry over $k$  does  $\gR$ remember? The answer seems to be far from being simple.

Let us also introduce another related ring, denoted $\gRs$, which is  generated by isomorphism classes of smooth complete $k$-varieties with product as multiplication and with relations $[\emptyset]=0$ and $[\B Y X]-[E]=[X]-[Y]$ where $\B Y X$ is the blowup of a smooth complete $k$-variety $X$ along smooth, complete $k$-variety $Y\subset X$ and $E$ is the exceptional divisor (after \cite{bittner01} these relations are  referred to as Bittner relations).

In \cite{bittner01} it was shown that:
\bthm\label{th:bittner}
Let $k$ be an algebraically closed field of characteristic zero. Then there is a natural isomorphism 
$$\gRs\mapright{\isom} \gR$$
given by sending $[X]\mapsto [X]$. 
\ethm

The proof of \Cref{th:bittner} given in \cite{bittner01} also works if $char(k)=p>0$ provided one knows that the \emph{weak factorization theorem} of  \cite{matsuki02} holds for $k$ and in particular, proof of that theorem requires that \emph{embedded resolution of singularities} (see \cite{hironaka64,abhyankar66})  holds for $k$ in all dimensions.

In characteristic zero it was shown in \citep{poonen02} that $\gR$ is not a domain. For this assertion over non-algebraically closed fields see \cite{kollar05} (I thank J\'{a}nos Koll\'{a}r for alerting me to this reference).   

As far as I am aware the question of whether (or not) the Grothendieck ring of varieties over a field of positive characteristics has been open for some time. The purpose of this note is to prove that $\gRs$ is not a domain (see \Cref{th:main}) if $p>13$ or $p=11$. 

My result (i.e. \Cref{th:main}) provides a canonical and even an elegant answer  to this question (and the answer is uniform in all the aforementioned characteristics).
My approach differs from \cite{poonen02} in the following way: a key result of \cite{larsen03} which is needed in \cite{poonen02} is not available in characteristic $p>0$; this difficulty is circumvented by working with ring $\gRs$  which is conjecturally isomorphic to $\gR$ (under availability of embedded resolution of singularities in all dimensions);  for $\gRs$  it is simpler to construct motivic measures (such as the Albanese measure); the other innovation is that I use a theorem of Pierre Deligne for constructing relations
in $\gRs$.

Let $\L=[\A^1]\in \gR$ be the isomorphism class of the affine line in $\gR$. In \citep{borisov2018} it was shown (assuming $k=\C$) that $\L$ is a zero divisor in $\gR$. In \Cref{asanuma-example}, I observe that there is an infinite collection of  natural candidate relations in $\gR$ which show that  $\L$ is a zero divisor (for any field $k$ with $char(k)=p>0$). But I do not know how to prove this at the moment (see \Cref{asanuma-example} for more precise statements).

In \Cref{th:pi1-is-a-motivic-measure} I  prove that the \'etale fundamental group of a smooth, proper variety (in any characteristic)  provides a motivic measure on the ring $\gRs$. In particular, the \'etale fundamental group is a motivic measure on the Grothendieck ring $\gR$ of varieties over complex numbers. One may think of this result as an algebraic manifestation of the anabelian nature of \'etale fundamental groups of (some) algebraic varieties.

It is a pleasure to thank Ravi Vakil for a stimulating series of lectures surrounding the Grothendieck ring $\gR$ at the Arizona Winter School (\href{http://swc.math.arizona.edu/aws/2015/index.html}{AWS 2015}) which served as a starting point for this note. It is also a pleasure to thank Research Institute of Mathematical Sciences (RIMS, Kyoto) and Shinichi Mochizuki  for hosting my visit to RIMS in Spring 2018 and for answering questions related to his construction of Albanese varieties.

\numberwithin{equation}{subsection}
\section{First theorem}
\subsection{Existence of zero divisors} 
The main theorem of this note is the following:
\bthm\label{th:main}
Let $k$ be an algebraically closed field of characteristic $p>0$.  
\label{th:main-1} If $p>13$ or $p=11$ then the isomorphism class of any supersingular elliptic curve is a zero divisor in $\gRs$. In particular,  $\gRs$ contains zero-divisors.
\ethm
\bp
Let $AV_k$ denote the multiplicative monoid of isomorphism classes of abelian varieties over $k$ and let $\Z[AV_k]$ denote the monoid ring.
To prove the theorem I construct a natural homomorphism (\emph{the Albanese motivic measure}) $$\alb:\gRs\to \Z[AV_k]$$ given by $[X]\to[\alb(X)]$. This construction is the content of \Cref{th:alb3} below. Now let me prove \Cref{th:main} assuming the construction of this homomorphism (see \Cref{th:alb3}).  

Recall from \cite[Chapter V, Theorem 4.1(c)]{silverman-arithmetic} that there are exactly 
\be 
\delta_p=\left[\frac{p-1}{12}\right]+\varepsilon_p
\ee isomorphism classes of supersingular elliptic curve over an algebraically closed field of characteristic $p>0$, see \citep[Chapter 5]{silverman-arithmetic} for the definition of $\varepsilon_p$ and at any rate note that $\delta_p\geq 2$ if and only if $p>13$ or $p=11$.

Suppose  that   $\delta_p\geq 2$. Let $E_1/k$ be any supersingular elliptic curve and let $E_2/k$ be any supersingular elliptic curve not isomorphic to $E_1$. Then by a beautiful theorem of Pierre Deligne (see \citep{oort87b}) there exists an isomorphism $E_1\times E_1\isom E_1\times E_2$. Thus, one has
\be 
[E_1\times E_1]=[E_1]\cdot [E_1]=[E_1]\cdot [E_2]=[E_1\times E_2].
\ee
in $\gRs$. In particular, we have 
\be 
[E_1]\cdot([E_1]-[E_2])=0.
\ee
Observe that as the images of  $[E_1],[E_1]-[E_2]$ are non-zero in $\Z[AV_k]$, and hence one sees that $[E_1],[E_1]-[E_2]\neq0\in\gRs$ so  both are zero divisors in $\gRs$. In particular, $[E_1]$ is a zero-divisor in $\gRs$ as claimed. This proves \Cref{th:main-1}.

\ep

\brem
More generally Deligne's Theorem (see \cite{oort87b}) asserts that if $n\geq2$ and $E_1,\ldots,E_{2n}$ are any supersingular elliptic curves then $E_1\times \cdots \times E_n\isom E_{n+1}\times \cdots \times E_{2n}$. So this theorem provides many more relations in $\gRs$. 
\erem

\brem 
Recently, using \Cref{th:main-1} and its proof, it has been shown in \cite{bot2025} that $\gRs$ also has nilpotent elements.
\erem

\subsection{Albanese variety functor as a motivic measure}
The following theorem, which constructs the Albanese motivic measure on $\gRs$ circumvents the deeper construction of motivic measures on $\gR$ due to \cite{larsen03} which is crucially dependent on (embedded) resolution of singularities.
\bthm\label{th:alb3} Let $k$ be an algebraically closed field of arbitrary characteristic. Let $AV_k$ be the multiplicative monoid of isomorphism classes of abelian $k$-varieties and let $\Z[AV_k]$ be the monoid ring. Then one has a natural morphism $$\alb{}:\gRs\to \Z[AV_k]$$ which is given by $[X]\to [\alb(X)]$. In particular, if $k$ has characteristic zero, then $\alb$ is a motivic measure on the Grothendieck ring of $\gR$ of $k$-varieties. 
\ethm

\bp 
Recall that if $X$ is a smooth, complete $k$-variety then $\alb(X)$ is an abelian variety equal to the dual of the reduced Picard scheme $\Pic(X)_{red}$ of $X$. For understanding properties of Albanese varieties of arbitrary varieties over arbitrary fields arbitrary characteristics readers may find \cite{ramachandran01}, \cite[Appendix, Page 66]{mochizuki-topics} useful. 

It is clear that $[X]\to[\alb(X)]\in {AV}_k$ can be extended linearly  (with $[\emptyset]\mapsto 0$) to define a homomorphism of groups from the free abelian group generated by isomorphism classes of smooth, complete $k$-varieties. By \Cref{le:alb2}, for complete $k$-varieties one has $[X]\times[Y]=[X\times_k Y]\mapsto[\alb(X\times_k Y)]=[\alb(X)]\times[\alb(Y)]$. So the mapping respects multiplication. By \Cref{le:alb1} $[\alb(X')]=[\alb(X)]$ and $[\alb(Y')]=[\alb(Y)]$. Thus the element $[X']-[X]-([E]-[Y])$ which maps to $[\alb(X')]-[\alb(X)]-([\alb(E)]-[\alb(Y)])=0$.  Thus this homomorphism  respects multiplication and Bittner relations and hence this homomorphism factors through the quotient ring $\gRs$. Hence one has the induced morphism $\gRs\to \Z[AV_k]$. The last assertion now follows from \Cref{th:bittner}.
\ep

\blem\label{le:alb1}
Suppose $k$ is algebraically closed and $Y\subset X$ be smooth, complete $k$-varieties. Let $X'=\B Y X$ and $Y'=E\subset X'$ be the exceptional divisor. Then one has  natural isomorphisms 
$$\alb(X')\isom \alb(X)$$ and $$\alb(Y')\isom\alb(Y).$$
\elem

\bp 
The proof follows from the well-known fact: if $A$ is an abelian variety then there are no non-constant morphisms from a projective space $\P^m\to A$ \cite{mumford-abelian2}. Suppose $A'=\alb(X')$ and $A=\alb(X)$. Since there are no non-constant morphisms from a projective space to an abelian variety, the tautological morphism $X'\to A'$ factors as $X'\to X\to A'$ and by the universal property of Albanese variety $A$ of $X$ this factors as $X'\to X\to A\to A'$ and in particular, it follows that any morphism from $X'$ to an abelian variety factors through $A$ so $A$ is the Albanese variety of $X'$ hence $A=A'$ and similarly one sees that $\alb(Y')=\alb(Y)$.
\ep

\blem\label{le:alb2}
Let $X_1,X_2$ be two smooth, complete $k$-varieties. Then $$\alb(X_1\times_k X_2)\isom \alb(X_1)\times \alb(X_2).$$
\elem

\bp 
This is immediate from the universal property of Albanese varieties.
\ep

\subsection{A relation which may be used to prove that $\L$ is a zero divisor in $\gR$}\label{asanuma-example}
By \cite{asanuma87,gupta2014} there exist a remarkable collection of explicitly defined  smooth affine algebras  $A/k$,
given explicitly by $A=k[x,y,z,w]/(x^my+z^{p^e}+w+w^{sp})$ where $e,m,s$ are positive integers such that $p^e\nmid sp$, $sp\nmid p^e$, and 
 with the following remarkable  properties: 
\beas
\spec(A)\times_k \A^1&\isom&\A^4\\
\spec(A)&\not\isom&\A^3.
\eeas
This example shows that in $\gR$ one has the relation
\be 
[\spec(A)]\cdot\L-\L^4=([\spec(A)]-\L^3)\cdot\L=0.
\ee 
and in particular, $\L$ is a zero-divisor if $$[\spec(A)]-\L^3\neq0.$$ 

But I do not know how to prove that $[\spec(A)]\neq\L^3$ in $\gR$. In fact suppose $A,A'$ are two such algebras with different values of parameters $e,m,s$ chosen so that $\spec(A),\spec(A')$ are not isomorphic as schemes then it is enough to prove that $[\spec(A)]\neq[\spec(A')]$ in $\gR$ as one also has the relation $$[\spec(A)]\cdot\L=\L^4=[\spec(A')]\cdot\L.$$ But again I do not know how to prove that $[\spec(A)]\neq[\spec(A')]$ in $\gR$ (or $\gRs$).

\newcommand{\prof}{\mathscr{Prof}}
\newcommand{\Bl}{{\rm Bl}}
\section{Second theorem: the \'etale fundamental group is a motivic measure}\numberwithin{equation}{section}
In this section, the base field $k$ is assumed to be algebraically closed of \emph{arbitrary characteristic}, and by the fundamental group of algebraic variety, I mean its \'etale fundamental group (see \cite[Expos{\'e} X]{sga1}). As is well-known, this is a profinite group (see \cite{sga1}). Let $X/k$ be a smooth, complete variety over a field. Then it is well-known that the isomorphism class of the \'etale fundamental group of $X/k$ is independent of the choice of the base point used to calculate the fundamental group and I will suppress the base point from my notation. In this note I record the following remark:
\bthm\label{th:pi1-is-a-motivic-measure}
Let $k$ be an algebraically closed field of arbitrary characteristic.  Let $\prof$ be multiplicative monoid of isomorphism classes of profinite groups. Then the mapping which maps a connected, smooth, proper variety  $X/k$ to its \'etale fundamental group $\pi_1(X/k)$ extends to a ring homomorphism:
$$\Pi:\gRs\to \Z[\prof].$$
In particular, if $k$ is of characteristic zero then $X\mapsto \pi_1(X/k)$ extends to a ring homomorphism 
$$\Pi: \gR \to \Z[\prof].$$ So the \'etale fundamental group is a motivic measure.
\ethm
\brem 
Since Grothendieck's anabelian program (and the results of A. Tamagawa, S. Mochizuki and others) suggests that the isomorphisms between \'etale fundamental groups should (in many reasonable cases of interest) rise from $\Z$-isomorphisms between the schemes,
one may (roughly) think of the homomorphism  $\Pi$, and especially its kernel, as an algebraic way of measuring anabelian properties of  $k$-varieties. We caution the reader that, in contrast to the anabelian program, in \Cref{th:pi1-is-a-motivic-measure}, $k$ is assumed to be algebraically closed. 
\erem
\bp 
For a profinite group $P$ let $[P]$ denote its isomorphism class. By \Cref{th:bittner}, it suffices to prove that if $X/k$ is a smooth, proper, connected variety and $Y\subset X$ is a smooth, subvariety of $X$,  and   $E\subset \Bl_Y(X)$ is the exceptional divisor of $\Bl_Y(X)\to X$ then it is sufficient to   prove that the following hold for the fundamental group: $$\pi_1(\Bl_Y(X)/k)\isom \pi_1(X/k),$$  and $\pi_1(E/k)\isom \pi_1(Y/k)$
and $$\pi_1(X_1\times_k X_2)\isom \pi_1(X_1/k)\times \pi_1(X_2/k).$$
For the first, one reduces to the case that $Y$ is connected as $\Bl_Y(X)$ can be obtained by successively blowing up connected components of $Y$. As the morphism $\Bl_Y(X)\to X$ is birational,  one can invoke \cite[Expos\'e X, Corollary 3.4]{sga1} to deduce $\pi_1(\Bl_Y(X)/k)\isom \pi_1(X/k)$. Note that $E\to Y$ is a projective bundle over $Y$ and as $\P^n$ is simply connected for all $n\geq 1$, one sees by \cite[Expos\'e X, Corollary 1.4]{sga1} that $\pi_1(E/k)\isom \pi_1(Y/k)$. The product relation follows from \cite[Expos\'e X, Corollary 1.7]{sga1}.   Hence one has the following equality of isomorphism classes of profinite groups:
\begin{align*}
{[\pi_1(\Bl_Y(X)/k)]}&= [\pi_1(X/k)], \\
[\pi_1(E/k)]&= [\pi_1(Y/k)],\\
[\pi_1(X_1\times_k X_2)] & =[\pi_1(X_1/k)]\times [\pi_1(X_2/k)].
\end{align*}
Hence one sees from the first two equations that the mapping $X\mapsto[\pi_1(X/k)]$ satisfies the Bittner relation
$$[\pi_1(\Bl_Y(X)/k)]-[\pi_1(E/k)]= [\pi_1(X/k)]-[\pi_1(Y/k)]$$
hold tautologically. This proves the theorem.
\ep
\bibliographystyle{plainnat}
\bibliography{grothendieck.bib, ../../master/master6.bib}
\end{document}